\documentclass{article}
\usepackage[utf8]{inputenc}

\usepackage{amssymb}
\usepackage{amsmath,mathtools}
\usepackage{amsthm}
\usepackage{amsfonts}
\usepackage{tikz}
\usepackage{xcolor}
\usepackage[top=1.0in, bottom=1.0in, left=1.0in, right=1.0in, headheight=1.0in]{geometry}
\usepackage{mathtools}
\usepackage{youngtab}
\usepackage{graphicx}
\usepackage{multicol}
\usepackage{paracol}
\graphicspath{ {./Graphs/} }
\usepackage[nottoc,notlof,notlot]{tocbibind}
\usepackage{changepage}

\newtheorem{theorem}{Theorem}[section]
\newtheorem{lemma}[theorem]{Lemma}

\newtheorem{prop}[theorem]{Proposition}

\newtheorem{corollary}[theorem]{Corollary}

\title{Graphs Identifiable by Degree Sequence and Chromatic Number}
\author{Rebecca Whitman\\
\small Department of Mathematics\\
\small University of California, Berkeley\\
\small Berkeley, CA 94720\\
\small\tt  rebecca\_whitman@berkeley.edu\
}
\date{\today}

\begin{document}

\maketitle

\begin{abstract}    
Unigraphs are graphs identifiable up to isomorphism from their degree sequences. Given a class $\mathcal{A}$ of graphs, we define the class of $\mathcal{A}$-unigraphs to be graphs identifiable from degree sequence and membership in $\mathcal{A}$. While these classes are often not hereditary, we provide characterizations of the largest hereditary subclass contained in the bipartite-unigraphs, the $k$-partite unigraphs, the perfect-unigraphs, and the chordal-unigraphs. We also characterize the largest hereditary subclass contained in the bipartite-unigraphs in terms of structure, degree sequence, and a partial order on degree sequences due to Rao. Lastly, we show that all unigraphs $G$ satisfy the bound $\chi(G) \le \omega(G) + 1$ and are hence apex-perfect graphs. 
\end{abstract}

\section{Introduction}

In this paper we consider a series of hereditary families related to unigraphs. All graphs are finite and simple. Given a graph $G = (V(G), E(G))$, the degree sequence $d$ of $G$ is the list of the degrees of the graph's vertices, written in non-increasing order. Where terms are repeated, we often denote their multiplicity with an exponent; for instance, $(3,2,2,2,2,1)$ and $(3, 2^4, 1)$ are the same degree sequence. We say that $G$ \emph{realizes} $d$ or is a \emph{realization} of $d$. Frequently, a degree sequence $d$ has numerous realizations, such that knowing $d$ is insufficient for determining much about a realization $G$. Where $d$ has exactly one realization up to isomorphism, we call $d$ \emph{unigraphic} and $G$ a \emph{unigraph}. 

Unigraphs have been studied for about fifty years, with structural and degree sequence characterizations by Tyshkevich and Chernyak (see \cite{TyCh78,TyCh79a, TyCh79b, TyCh79c} and the English translation \cite{Ty00}). The class of \emph{split graphs} - graphs for which the vertex set $V(G)$ can be partitioned into a subset $K$ inducing a clique and a subset $S$ inducing a stable set - is an important tool in these characterizations of the class of unigraphs. Specifically, both characterizations rely on an operation defined by Tyshkevich that decomposes a graph into $n$ split graph components and one unspecified component \cite{Ty80} \cite{Ty00}. The authors are then able to classify exactly the indecomposable components of a unigraph by structure and degree sequence. We summarize results pertaining to unigraphs here, and defer a fuller treatment of the decomposition operation to Section \ref{sec:coloring}.

\begin{theorem} \cite{Ty00}
\label{prop:Tyshkevich_unigraph}
    A graph $G$ is a unigraph if and only if each component of its decomposition is a unigraph. 
\end{theorem}

Given a graph $G$, we denote its complement by $\overline{G}$. Given a split graph $G$, a \emph{$KS$-partition} is a partition of $V(G)$ into subsets $K$ inducing a clique and $S$ inducing a stable set. We write $(G,K,S)$ to indicate a split graph $G$ with $KS$-partition $K \cup S$. Given $(G,K,S)$, the \emph{inverse} $G^I$ is the split graph obtained from $G$ by removing all edges with both endpoints in $K$ and adding all edges with both endpoints in $S$ (thus switching the clique and stable set of the graph). Note that the inverse depends on the choice of $KS$-partition. 

Neither complementation nor inversion changes a graph's status as a unigraph. As a result, indecomposable unigraphs can be characterized as follows.

\begin{theorem} \cite{Ty00}
\label{prop:indecomposable_unigraphs}
    A graph $G$ is an indecomposable unigraph if and only if $G$, $\overline{G}$, or, if split, $G^I$ or $\overline{G^I}$ are isomorphic to $K_1$ or $C_5$ or belong to one of seven infinite families of graphs. 
\end{theorem}

Included in these seven families of unigraphs are all graphs isomorphic to the union of a star and $nK_2$, $n \ge 0$, and split graphs where all vertices in the stable set are of degree $1$. We refer the reader to \cite{Ty00} for the complete list.

A class $\mathcal{A}$ of graphs is \emph{hereditary} if it is closed under taking induced subgraphs of its members. All hereditary classes can be characterized by a set of forbidden induced subgraphs. For any class $\mathcal{A}$ of graphs, we define two closely related hereditary classes. Let $H(\mathcal{A})$ be the \emph{hereditary subclass of $\mathcal{A}$}, the largest hereditary class contained within $\mathcal{A}$. For example, the perfect graphs are the hereditary subclass of the set of graphs with equal clique and chromatic numbers. The set $H(\mathcal{A})$ is equivalently the subset of elements of $\mathcal{A}$ for which all induced subgraphs are also elements of $\mathcal{A}$. We call a graph $G$ in $H(\mathcal{A})$ a \emph{hereditary $\mathcal{A}$-graph} (e.g., hereditary unigraph). Let $HC(\mathcal{A})$ be the \emph{hereditary closure of $\mathcal{A}$}, the smallest hereditary class containing $\mathcal{A}$. For example, the class of forests is the hereditary closure of the class of trees. The set $HC(\mathcal{A})$ is equivalently the set of all graphs induced in an element of $\mathcal{A}$. 

\begin{figure}
\centering
  \includegraphics[height=4cm]{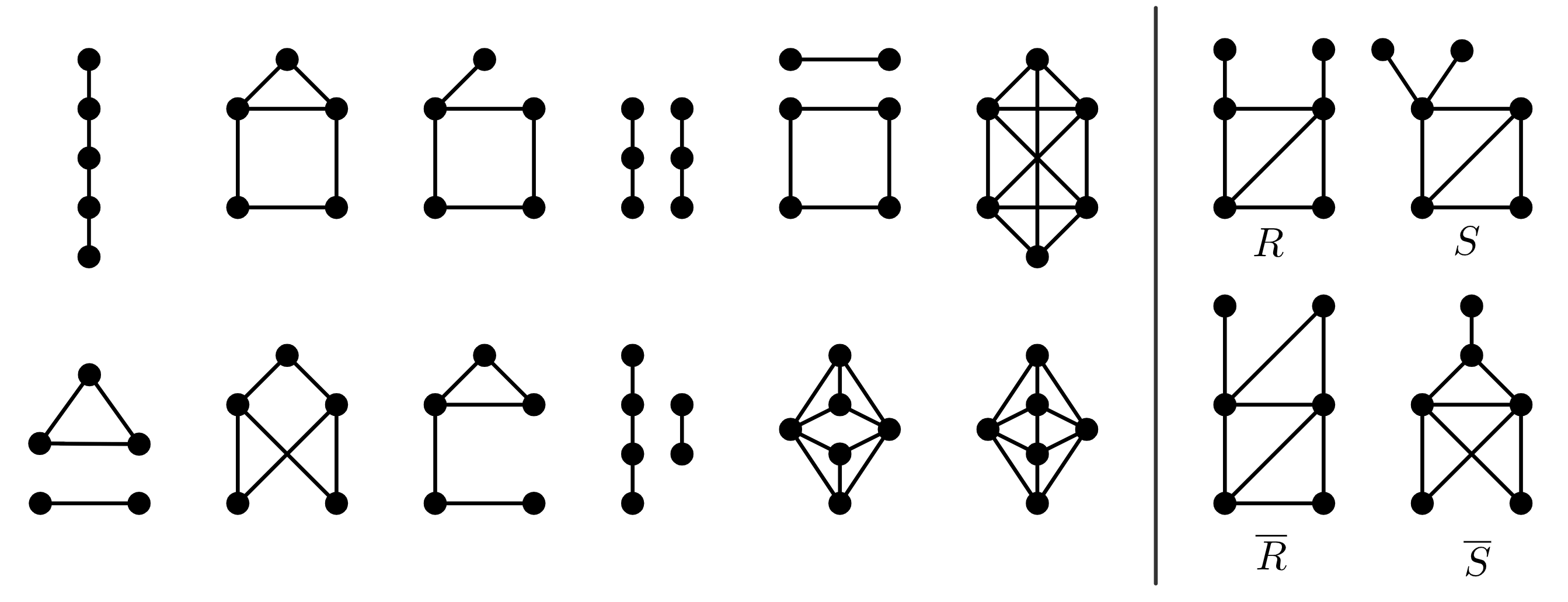}
  \caption{The 16 forbidden induced subgraphs of the hereditary unigraphs.}
  \label{fig:hereditary_unigraphs}
\end{figure}

The class $\mathcal{U}$ of unigraphs contains a number of important hereditary classes, including threshold \cite{HaSi81}, matroidal \cite{MaEtAl84}, and matrogenic graphs \cite{Ty84},
but is not itself hereditary. The smallest counterexample is the unigraphic degree sequence $(3,3,3,3,3,1)$: its realization contains two non-isomorphic realizations of $(3,2,2,2,1)$. Barrus \cite{Ba12} \cite{Ba13} studied the hereditary subclass of the unigraphs and characterized them in terms of their degree sequences, their structure, and a finite list of forbidden induced subgraphs for the class. Specifically, there are sixteen forbidden induced subgraphs, each containing five to seven vertices, which are shown in Figure \ref{fig:hereditary_unigraphs}. Let $\mathcal{F}_{H(\mathcal{U})}$ be the set of these graphs. The graphs $R$, $S$, and their complements, defined as shown in Figure \ref{fig:hereditary_unigraphs} are the four split graphs in this set. 

\begin{theorem}\cite{Ba12}
\label{prop:hereditary_unigraphs}
    A graph $G$ is a hereditary unigraph if and only if $G$ contains no element of $\mathcal{F}_{H(\mathcal{U})}$ as an induced subgraph. 
\end{theorem}

Barrus, Trenk, and Whitman \cite{BaTrWh24} similarly characterized the hereditary closure of the unigraphs, in terms of degree sequences, structure, and a finite list of forbidden induced subgraphs. They also provide a characterization of $H(\mathcal{U})$ and $HC(\mathcal{U})$ in terms of minimally forbidden degree sequences, in accordance with Rao \cite{Rao80}. 

In this paper, we identify unique realizations within a number of important classes of graphs. In other words, what graphs - within the restricted universe of a particular graph class - are the only realizations of their degree sequences? Define an \emph{$\mathcal{A}$-unigraph} to be a graph $G \in \mathcal{A}$ such that its degree sequence $d$ has exactly one realization that is an element of $\mathcal{A}$. We call $d$ \emph{$\mathcal{A}$-unigraphic}. Note that $\mathcal{A}$-unigraphs and $\mathcal{A}$-unigraphic degree sequences are not necessarily unigraphs and unigraphic degree sequences without qualification. For example, $P_5$ is a bipartite-unigraph, since the only other realization of $d(P_5) = (2,2,2,1,1)$ is $K_3 + K_2$, which is not a bipartite graph, but $P_5$ is not a unigraph writ large. 

\begin{figure}
\centering
  \includegraphics[height=3cm]{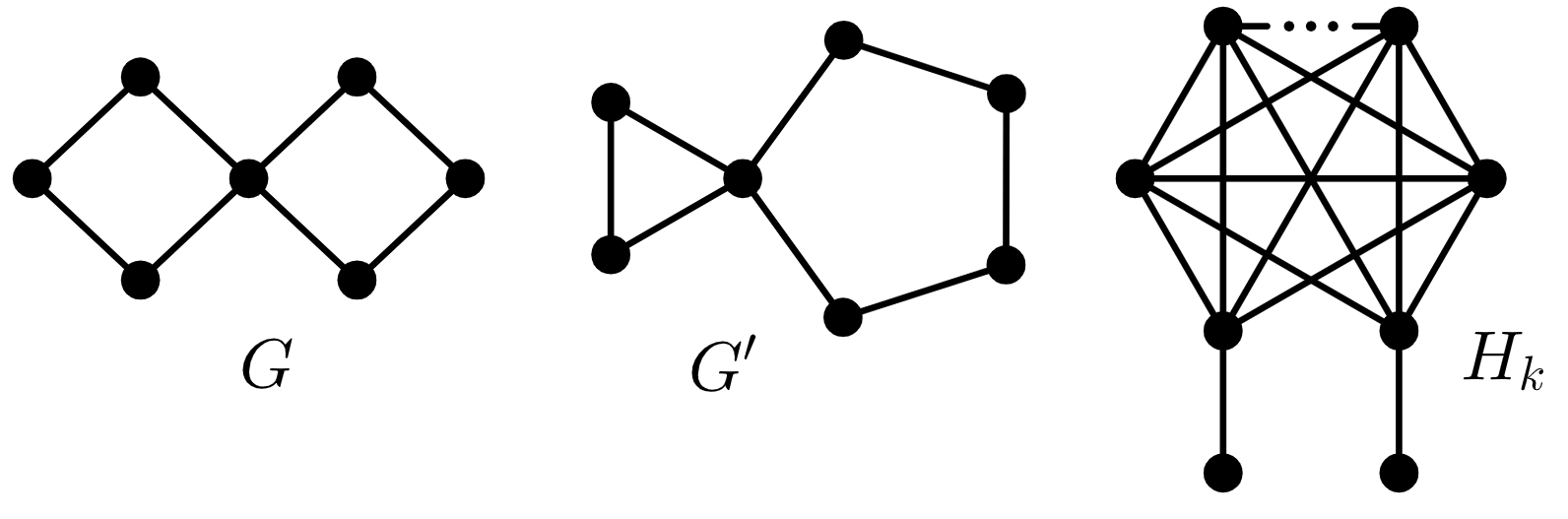}
  \caption{A family of $k$-partite graphs that are not $k$-partite-unigraphs.}
  \label{fig:counterexample}
\end{figure}

We generalize bipartite-unigraphs to the families $\mathcal{U}_k$ of $k$-partite-unigraphs for $k \ge 2$. Each is the set of graphs exactly determinable from degree sequence and an upper bound $k$ on the chromatic number $\chi$ of the graph. For example, $P_5$ is not a $k$-partite-unigraph for $k \ge 3$, since $K_3+K_2$ is also $k$-partite. While $k$-partite graphs form a hereditary class, $\mathcal{U}_k$ does not. For $k=2$, consider the two realizations $G$ and $G'$ of degree sequence $(4, 2^6)$, both shown in Figure \ref{fig:counterexample}. Since $G'$ contains $K_3$, $G$ is a bipartite-unigraph. However, $G$ contains $2P_3$, which is not a bipartite-unigraph (the other realization of its degree sequence is $P_4+K_2$, which is also bipartite), so $G$ is not a hereditary bipartite-unigraph. For $k \ge 3$, consider the degree sequence $d_k = (k^{k+1},1,1)$. The sequence has exactly two realizations: one, the graph $H_k$ shown in Figure \ref{fig:counterexample}, is $k$-partite, and the other, $K_{k+1} + K_2$, is not. However, $H_k$ contains $P_5$ as an induced subgraph, so $H_k$ is not a hereditary $k$-partite-unigraph. 

We provide a partial forbidden induced subgraph characterization of $H(\mathcal{U}_k)$, the hereditary subclass of the $k$-partite-unigraphs. For $H(\mathcal{U}_2)$, we give exact structural, degree sequence, and forbidden induced subgraph characterizations. In particular, we show that hereditary bipartite-unigraphs are closely related to complete bipartite graphs. We also characterize a set of forbidden partitioned degree sequences for $H(\mathcal{U}_2)$, \`{a} la Rao \cite{Rao80}.

From here we turn to a more general discussion of coloring unigraphs. In fact, unigraphs also have excellent graph coloring properties. Let $\omega(G)$ and $\chi(G)$ denote, respectively, the clique number and chromatic number of a graph $G$. For a hereditary class $\mathcal{A}$, $f: \mathbb{N} \rightarrow \mathbb{N}$ is a $\chi$-bounding function on $\mathcal{A}$ if for all $G \in \mathcal{A}$, $\chi(G) \le f(\omega(G))$. For perfect graphs, the function $f(x) = x$ is $\chi$-bounding, and this is of course best possible. See \cite{ScSe20} for a survey of hereditary classes with known $\chi$-bounding functions. We show that the class $\mathcal{U}$ of unigraphs is bounded by $f(x) = x+1$, pointing to a strong relationship between graph realization and graph coloring. 

In the final section of the paper, we expand on this relationship to consider two other families of hereditary $\mathcal{A}$-unigraphs, again with stronger coloring properties: the hereditary perfect-unigraphs $H(\mathcal{U}_P)$ and the hereditary chordal-unigraphs $H(\mathcal{U}_C)$. We provide a forbidden induced subgraph characterization of each class. 
\section{Hereditary $k$-partite-unigraphs}

We begin with four equivalent characterizations of the hereditary subclass of bipartite graphs, $H(\mathcal{U}_2)$: by degree sequence, by graph structure, by forbidden induced subgraphs, and by Rao-forbidden partitioned degree sequences. 

Before presenting these characterizations in Theorem \ref{prop:hbu_characterization}, we construct a partial ordering of degree sequences of bipartite graphs. 

An unordered pair $\{d_1,d_2\}$ is a \emph{bipartitioned degree sequence pair} if there exists a bipartite graph $G$ with bipartition $V(G) = A \cup B$ such that $d(A) = d_1$ and $d(B) = d_2$. Thus $d(G) = d_1 \cup d_2$. 

Given bipartitioned pairs $\{d_1, d_2\}$ and $\{e_1, e_2\}$, we say that $\{d_1, d_2\}$ \emph{bipartite-Rao-contains} $\{e_1, e_2\}$ (and write $\{d_1, d_2\} \succeq \{e_1, e_2\}$) if and only if there exist bipartite graphs $G$ and $H$ such that (i) $V(G) = G_1 \cup G_2$ and $V(H) = H_1 \cup H_2$ are bipartitions of $G$ and $H$, respectively, with $d(G_1) = d_1$, $d(G_2) = d_2$, $d(H_1) = e_1$, and $d(H_2) = e_2$, and (ii) there exists an injective homomorphism $\phi: V(H) \rightarrow V(G)$ with $\phi(h) \in G_1$ if and only if $h \in H_1$. We also say that a degree sequence $d$ bipartite-Rao-contains a bipartitioned pair $\{e_1,e_2\}$ if there exists a bipartition $\{d_1,d_2\}$ of $d$ with $\{d_1, d_2\} \succeq \{e_1, e_2\}$. 

A bipartitioned pair is a \emph{forbidden pair} if any bipartite realization of the sequence is not in $H(\mathcal{B}_2)$. Among forbidden pairs, a  forbidden pair is \emph{bipartite-Rao-minimal} if it does not bipartite-Rao-contain any other forbidden pair. Forbidden pairs and forbidden induced subgraphs are related per the following lemma. 

\begin{lemma}
\label{prop:Rao_forbidden_containment}
    A forbidden pair is bipartite-Rao-minimal if and only if all its realizations are forbidden induced subgraphs. 
\end{lemma}
\begin{proof}
    Suppose that $\{d_1, d_2\}$ is a forbidden pair. If $\{d_1, d_2\}$ is not bipartite-Rao-minimal, then there exists a forbidden pair $\{e_1, e_2\}$ bipartite-Rao-contained in $\{d_1, d_2\}$. Thus, there exist bipartite realizations $E, E'$ of $\{e_1, e_2\}$ and $D$ of $\{d_1, d_2\}$ with $E$ induced in $D$, and we conclude that $D$ is not a forbidden induced subgraph. 

    If instead $\{d_1, d_2\}$ has a realization $D$ that is not a forbidden induced subgraph, it follows that $D$ contains some forbidden induced subgraph $E$. The forbidden pair corresponding to $E$ is thus bipartite-Rao-contained in $\{d_1, d_2\}$, which we conclude is not bipartite-Rao-minimal. 
\end{proof}

\begin{figure}
\centering
  \includegraphics[height=4.5cm]{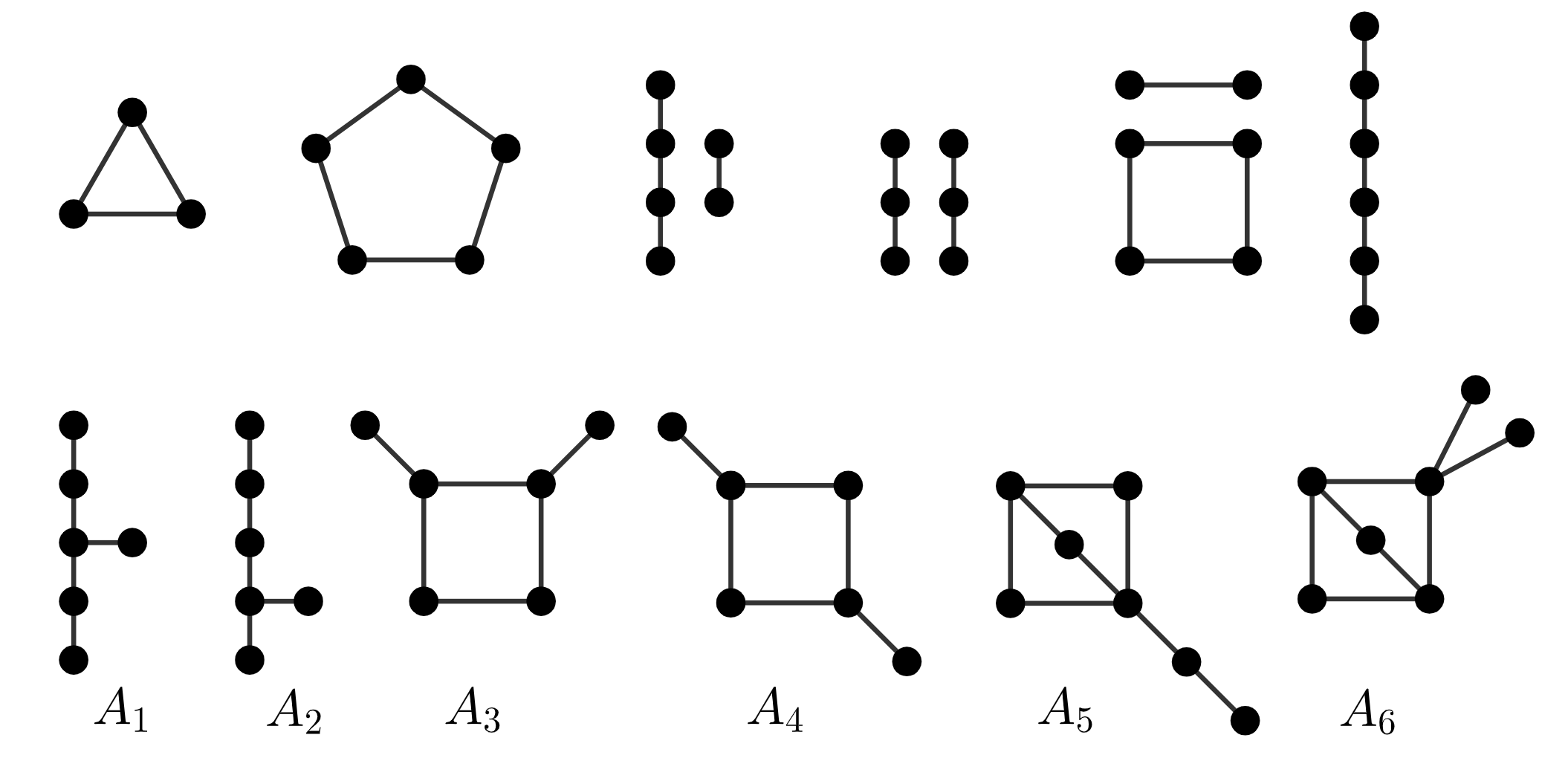}
  \caption{The class of hereditary bipartite-unigraphs has twelve forbidden induced subgraphs.}
  \label{fig:bipartite_mfis}
\end{figure}

We will thus be able to characterize $H(\mathcal{U}_2)$ by its forbidden pairs. Let $\mathcal{F}_{H(\mathcal{U}_2)} = \{K_3, C_5, P_4+K_2, 2P_3, C_4 + K_2, P_6, A_1, A_2, A_3, A_4, A_5, A_6 \}$, as shown in Figure \ref{fig:bipartite_mfis}. Let $\mathcal{R}_{H(\mathcal{U}_2)}$$ = \{ \{(2,2),(1,1,1,1)\}$, \\
        $\{(2,1,1),(2,1,1)\}$,
        $\{(2,2,1),(2,2,1)\}$,
        $\{(3,1,1),(2,2,1)\}$,
        $\{(3,2,1),(3,2,1)\}$,
        $\{(4,3,1),(2,2,2,2)\}$, \\
        $\{(4,2,2),(3,3,1,1)\} \}$. 
We show, among other characterizations, that $\mathcal{F}_{H(\mathcal{U}_2)}$ is the set of forbidden induced subgraphs of ${H(\mathcal{U}_2)}$, and that $\mathcal{R}_{H(\mathcal{U}_2)}$ is the set of bipartite-Rao-minimal pairs. 

Note that a graph $G$ is a hereditary bipartite-unigraph if and only if for any $v \not \in V(G)$, $G \cup \{v\}$ is a hereditary bipartite-unigraph. Hence we assume in the following characterization that $G$ contains no isolated vertices. 

\begin{theorem}
\label{prop:hbu_characterization}
    The following are equivalent for a graph $G$ containing no isolated vertices: 
    \begin{enumerate}
        \item [(i)] $G$ is a hereditary bipartite-unigraph. 
        \item [(ii)] $G$ contains no element of $\mathcal{F}_{H(\mathcal{U}_2)}$ as an induced subgraph.
        \item [(iii)] $G$ is bipartite and there exists a vertex $v \in V(G)$ such that either $G$ or $G - \{v\}$ is isomorphic to $K_{a,b}$ for some $a,b \ge 0$. 
        \item [(iv)] $G$ is bipartite and its degree sequence $d$ is equal to $(a^c, (a-1)^{b-c},b^{a-1},c)$ for some $a,b,c \ge 0$. 
        \item [(v)] $G$ is bipartite and its degree sequence $d$ does not bipartite-Rao-contain any element of $\mathcal{R}_{H(\mathcal{U}_2)}$. 
        \end{enumerate}
\end{theorem}
\begin{proof}
    (i) $\rightarrow$ (ii): 
    
    Suppose $G \in $ $H(\mathcal{U}_2)$. Since $G$ is bipartite, it contains no $K_3$ or $C_5$. The pairs of graphs $P_4+K_2$ and $2P_3$, $C_4+K_2$ and $P_6$, $A_1$ and $A_2$, and $A_3$ and $A_4$ realize the same degree sequences, so are not themselves bipartite-unigraphs and cannot be induced subgraphs of $G$. The degree sequences $(4,3,2,2,2,2,1)$ of $A_5$ and $(4,3,3,2,2,1,1)$ of $A_6$ each also admit multiple bipartite realizations, so neither $A_5$ nor $A_6$ is contained in $G$. 

    (ii) $\rightarrow$ (iii): 
    
    Suppose $G$ contains no element of $\mathcal{F}_{H(\mathcal{U}_2)}$ as an induced subgraph. Note that all indecomposable unigraphs which are bipartite, as identified in Theorem \ref{prop:indecomposable_unigraphs}, satisfy (iii).

    If $G$ is disconnected, then since $2P_3$ and $K_3$ are forbidden, at most one component, which we call $H$, has $3$ or more vertices. All other components are isomorphic to $K_2$, so $H$ cannot contain $P_4$ or $C_4$ as an induced subgraph. If $H$ has only two vertices, then $G$ is isomorphic to $nK_2$ and is a unigraph by \ref{prop:indecomposable_unigraphs}. If $H$ has at least three vertices, then $H$ contains an induced $P_3$, on, say, vertex set $\{a,b,c\}$ with edges $ab$ and $bc$. All other vertices in $H$ must have neighbor set exactly $\{b\}$, else $P_4$, $C_4$, or $K_3$ is produced as an induced subgraph. Thus $H$ is isomorphic to a star, and $G$, comprising a star and $nK_2$, is a unigraph by \ref{prop:indecomposable_unigraphs}.

    Otherwise, suppose $G$ is connected, and has at least $8$ vertices (it is trivial to check that all graphs on $7$ or fewer vertices contain no induced subgraph from $\mathcal{F}_{H(\mathcal{U}_2)}$). If $G$ contains no $P_3$ as an induced subgraph, then $G$ is isomorphic to $K_2$, for a contradiction. 
    
    If the longest path in $G$ has three vertices, let $\{a,b,c\}$ with edges $ab$ and $bc$ induce such a path. Either all additional vertices in $G$ have neighbor set $\{b\}$ and $G$ is a star and hence a unigraph by \ref{prop:indecomposable_unigraphs}, or $G$ contains some vertex adjacent to $a$ and $c$ but not $b$. Let $AC$ be the set of vertices in $V(G) - \{b\}$ adjacent to $a$ and $c$ but not $b$, and $B$ be the set of vertices in $V(G) - \{a,c\}$ adjacent to $b$ but not $a$ or $c$. If $AC$ and $B$ are complete to one another or missing at most one edge, then $G$ meets condition (iii). If $|AC| \le 1$ and $B$ contains at least two vertices not adjacent to a vertex in $AC$, then $G$ meets condition (iii). If $|AC| \ge 2$ and $B$ contains at least two different vertices not adjacent to the same vertex in $AC$, then $G$ contains an induced $A_6$. If $|AC| \ge 2$ and $B$ contains at least two different vertices not adjacent to different vertices in $AC$, then $G$ contains an induced $A_4$. If $|AC| \ge 2$ and $B$ contains only one vertex not complete to $AC$, but said vertex is non-adjacent to at least two vertices in $AC$, then $G$ meets condition (iii). Hence if the longest path in $G$ has three vertices, $G$ is a bipartite-unigraph. 

    If the longest path in $G$ has four vertices, let $\{a,b,c,d\}$ induce such a path with edges $ab, bc,$ and $cd$. By excluding $P_5$, $K_3$, and $C_5$, $G$ may contain a set $B$ of vertices adjacent to $b$ but not $a$, $c$, or $d$; a set $C$ of vertices adjacent to $c$ but not $a$, $b$, or $d$; a set $AC$ of vertices adjacent to $a$ and $c$ but not $b$ or $d$; and a set $BD$ of vertices adjacent to $b$ and $d$ but not $a$ or $c$. Since $G$ contains no induced $A_3$, $B$ is complete to $AC$ and anti-complete to $C$, and $C$ is complete to $BD$. If $|AC| \ge 2$, it is complete to $BD$ and if $|BD| \ge 2$, it is complete to $AC$, else $G$ contains an induced $A_2$ or $A_4$. If $B$ and $C$ are non-empty, then $AC$ and $BD$ are both empty lest $G$ contain $A_3$. Now, however, $G$ is a unigraph by Theorem \ref{prop:indecomposable_unigraphs}. Otherwise assume without loss of generality that $C$ is empty. Since $B \cup \{a, c\} \cup BD$ is complete to $AC \cup \{b\}$, $G$ meets condition (iii).

    If the longest path in $G$ has five vertices, let $\{a,b,c,d,e\}$ induce such a path with edges $ab, bc, cd$ and $de$. Since $G$ contains no $K_3, C_5, P_6, A_1, A_2, A_3,$ or $A_4$, $G$ may contain a set $AC$ of vertices adjacent to $a$ and $c$ but not $b$, $d$, or $e$; a set $CE$ of vertices adjacent to $c$ and $e$ but not $a$, $b$, or $d$; a set $AE$ of vertices adjacent to $a$ and $e$ but not $b$, $c$, or $d$; and a set $ACE$ of vertices adjacent to $a$, $c$, and $e$ but not $b$ or $d$. The set $AC \cup CE \cup AE \cup ACE$ is a stable set. If $|AC| \ge 2$, then $G$ contains an induced $A_5$. If $|CE| \ge 2$, then $G$ contains an induced $A_5$. If $|AE| \ge 2$, then $G$ contains an induced $A_4$. If $|ACE| \ge 2$, then $G$ contains an induced $A_4$. If $AC$ and $AE$ are both non-empty, or $CE$ and $AE$ are both non-empty, then $G$ contains an induced $A_2$. If $AE$ and $ACE$ are both non-empty, then $G$ contains an induced $A_4$. If $AC, CE$, and $ACE$ are all non-empty, then $G$ contains an induced $A_2$. Hence $|G| \le 8$, for a contradiction. 

    (iii) $\rightarrow$ (iv): 

    Assume $G$ has no isolated vertices. If $G$ is isomorphic to $K_{a,b}$, then its degree sequence is $(a^b, b^a)$, meeting condition (iv) with $c = b$. Otherwise, there exists $v \in V(G)$ such that $G - \{v\}$ is isomorphic to $K_{a,b}$ for some $a,b \ge 0$. Let $G - \{v\}$ have bipartition $V(G) - \{v\} = A \cup B$, where the vertices of $A$ are of degree $b$ and the vertices of $B$ are of degree $a$. Let $\deg(v) = x$. Since $G$ is bipartite, the neighbor set of $v$ is a subset of either $A$ or $B$; we assume $B$, without loss of generality. Thus the degree sequence of $G$ is $((a+1)^x, a^{b-x}, b^a, x)$, meeting condition (iv) (substituting $a+1$ for $a$). 

    (iv) $\rightarrow$ (i): 
    
    Let $G$ be bipartite and have degree sequence $d = (a^c, (a-1)^{b-c},b^{a-1},c,0^e)$. If $c=0$, then $d = (a^b, b^a)$ for some $a \ge b \ge 0$. Let $v_1, \ldots, v_b$ have degree $a$, and $v_{b+1}, \ldots, v_{a+b}$ have degree $b$. Up to symmetry, suppose $v_1 \in X$. Thus there are at least $a$ vertices in $Y$; since $a \ge b$, some vertex in $Y$ must have degree $b$. Thus there are at least $b$ vertices in $X$. Hence $|X| = b$ and $|Y| = a$. If $a > b$, it follows that $X = \{v_1, \ldots, v_b\}$ and $Y = \{v_{b+1}, \ldots, v_{a+b}\}$, and if $a = b$, then both $X$ and $Y$ contain $a$ vertices of degree $a$. In either case, the bipartition is unique, and all edges of the graph are fixed. 
    
    If $G$ is a complete bipartite graph with one additional vertex $v_1$ of degree $1 \le x \le a$, then either $d = (a^x, (a-1)^{b-x},b^{a-1},x)$ or $a>b$ and $d = (a^{b-1},x,b^x,(b-1)^{a-x})$. 

    In the first case, let a vertex $v_2$ of degree $a$ be in $X$. Thus there are at least $a$ vertices in $Y$, so some vertex in $Y$ has degree at least $b$. (This assumes $a \ge 2$; if $a=1$, then $|V(G)| \le 2$ and the proof is trivial.) Hence $|X| = b$ and $|Y| = a$. If $a=b$, then up to symmetry $v_1 \in X$, so $X$ contains $v_1$ and $a-1$ vertices of degree $a$, and $Y$ contains $x$ vertices of degree $a$ and $a-x$ vertices of degree $a-1$. Since all vertices of degree $a$ are complete to the opposite partition, there is exactly one way to realize the sequence as a bipartite graph. 

    If $a > b$, then all $x$ vertices of degree $a$ are in $X$. Since $\sum d = 2ab+2x-2b$, it follows that $\sum _{v \in X} \deg(v) = ab+x-b$, so the remaining $b-x$ vertices of $X$ must have total degree $ab-ax+x-b$. We conclude each has degree $a-1$, so $Y$ then contains $a-1$ vertices of degree $b$ and one vertex, $v_1$, of degree $x$. The vertices of degree $b$ in $Y$ are complete to $X$, so $v_1$ must be adjacent exactly to the vertices of degree $a$, such that there is exactly one bipartite realization. 

    In the second case, we assume $a>b$ and that some vertex $v_2$ of degree $a$ is in $X$. With at least $a$ vertices in $Y$, one must have degree at least $b$, so $|X| \ge b$. This implies $|X| = b$ and $|Y| = a$. Thus all $b-1$ vertices of degree $a$ are in $X$. Since $\sum d = 2ab+2x-2a$, it follows that $\sum _{v \in X} \deg(v) = ab+x-a$, so the remaining vertex of $X$ has degree $x$. Thus $Y$ contains exactly the $x$ vertices of degree $b$ and the $a-x$ vertices of degree $a$. As above, there is exactly one realization of this bipartitioned degree sequence. 

    If $G$ is a bipartite graph which is complete together with up to one additional vertex, then the same is true of all induced subgraphs $H$ of $G$. By the above argument all such $H$ are bipartite-unigraphs, and we conclude that $G$ is a hereditary bipartite-unigraph. 

    (v) $\rightarrow$ (iii): 
    
    Suppose that $G$ contains a graph from $\mathcal{F}_{H(\mathcal{U}_2)}$ as an induced subgraph. If $G$ contains $K_3$ or $C_5$, then it is not bipartite, for a contradiction. Otherwise, if $G$ contains $P_4+K_2$ or $2P_3$, then $d(G)$ bipartite-Rao-contains $\{(2,2),(1,1,1,1)\}$ and $\{(2,1,1),(2,1,1)\}$. If $G$ contains $C_4+K_2$ or $P_6$, then $d(G)$ bipartite-Rao-contains $\{(2,2,1),(2,2,1)\}$. If $G$ contains $A_1$ or $A_2$, then $d(G)$ bipartite-Rao-contains $\{(3,1,1),(2,2,1)\}$. If $G$ contains $A_3$ or $A_4$, then $d(G)$ bipartite-Rao-contains $\{(3,2,1),(3,2,1)\}$. If $G$ contains $A_5$, then $d(G)$ bipartite-Rao-contains $\{(4,3,1),(2,2,2,2)\}$. If $G$ contains $A_6$, then $d(G)$ bipartite-Rao-contains $\{(4,2,2),(3,3,1,1)\}$.  

    (iii) $\rightarrow$ (v):
    
    Suppose $G$ is bipartite and its degree sequence $d(G)$ bipartite-Rao-contains a sequence in $\mathcal{R}_{H(\mathcal{U}_2)}$. Thus there exists a bipartite graph $H$ realizing a bipartitioned degree sequence in $\mathcal{R}_{H(\mathcal{U}_2)}$, with $H$ induced in a bipartite realization $G'$ of $d(G)$. By Lemma \ref{prop:Rao_forbidden_containment}, all realizations of a sequence in $\mathcal{R}_{H(\mathcal{U}_2)}$ are forbidden induced subgraphs, which we showed above are precisely the graphs in $\mathcal{F}_{H(\mathcal{U}_2)}$. Hence $G'$ contains a graph from $\mathcal{F}_{H(\mathcal{U}_2)}$ as an induced subgraph. Every realization of the degree sequence of a graph in $\mathcal{F}_{H(\mathcal{U}_2)}$ itself contains an induced subgraph from $\mathcal{F}_{H(\mathcal{U}_2)}$. This implies $G$ contains an element of $\mathcal{F}_{H(\mathcal{U}_2)}$ as an induced subgraph, completing the proof. 
\end{proof}

Given a hereditary class $\mathcal{A}$ of graphs, the class of apex-$\mathcal{A}$ graphs is the set of graphs $G$ that induce an element of $\mathcal{A}$ on all but one vertex. Characterization (iii) is equivalent to stating that hereditary bipartite-unigraphs are apex-complete-bipartite.

We now turn to the study of $H(\mathcal{U}_k)$ for $k \ge 3$. These classes are not as interesting as $H(\mathcal{U}_2)$, and in fact hereditary $k$-partite-unigraphs are exactly hereditary unigraphs that are also $k$-partite. 

\begin{theorem}
\label{prop:multipartite}
    A graph $G$ is a hereditary $k$-partite-unigraph for $k \ge 3$ if and only if $G$ is $k$-partite and a hereditary unigraph.
\end{theorem}
\begin{proof}
    Suppose $G$ is a hereditary $k$-partite-unigraph. Of course $G$ is $k$-partite. The graphs of $\mathcal{L}$ are each $k$-partite for $k \ge 3$. Furthermore, they and all other realizations of their degree sequences are non-unigraphs by Theorem \ref{prop:hereditary_unigraphs}, so $G$ cannot contain any as an induced subgraph, so $G$ is a hereditary unigraph. We note the same does not hold for bipartite-unigraphs.  

    If $G$ is not a hereditary $k$-partite-unigraph, then either $G$ is not $k$-partite, or $G$ is $k$-partite but not a hereditary $k$-partite-unigraph.    
    Thus some induced subgraph $H$ of $G$ is $k$-partite but not a $k$-partite-unigraph. This implies $H$ is not a unigraph, so $G$ is not a hereditary unigraph.
\end{proof}

Let $\mathcal{F}_k$ be the set of forbidden induced subgraphs of the set of $k$-partite graphs. No known characterization exists for $\mathcal{F}_k$, but its smallest element is $K_{k+1}$. For $H(\mathcal{U}_k)$, we give a partial forbidden induced subgraph characterization, dependent on knowing $\mathcal{F}_k$. 

\begin{corollary}
    A graph $G$ is a hereditary $k$-partite-unigraph for $k \ge 3$ if and only if $G$ contains neither a graph from $\mathcal{F}_k$ nor a graph from $\mathcal{F}_{H(\mathcal{U})}$ as an induced subgraph. 
\end{corollary}

\section{Coloring Hereditary Unigraphs}
\label{sec:coloring}

Given a graph $G_0$ and a split graph $G_1$ with a partition of $V(G_1)$ into a clique $K$ and a stable set $S$, the \emph{composition of $(G_1, K, S)$ and $G_0$}, denoted $(G_1, K, S) \circ G_0$, is the graph  $G$ with vertex set $V(G_1) \cup V(G_0)$ and edge set $E(G_1) \cup E(G_0) \cup \{kv| k \in K, v \in V(G_0)\}$. Here $G$ is the union of the original graphs together with all edges from $G_1$'s clique to $G_0$. Note that the composition may vary depending on the choice of partition of $V(G_1)$ into a clique and stable set. Where $G_1$ admits exactly partition into a clique and a stable set, we write $G_1 \cup G_0$ for convenience. 

Composition is associative, so the multipart composition $(G_n, K^n, S^n) \circ \ldots \circ (G_1, K^1, S^1) \circ G_0$ is well-defined where $G_n, \ldots, G_1$ are split graphs with the given $KS$-partitions. 

Clique numbers and chromatic numbers behave predictably with respect to graph composition. 

\begin{prop}
\label{prop:composition_coloring}
    Let $G_i$, $1 \le i \le n$ be a split graph with fixed $KS$-partition $K^i \cup S^i$ and $G_0$ be a nonempty graph. Given $G = (G_n, K^n, S^n) \circ \ldots \circ (G_1, K^1, S^1) \circ G_0$, it holds that $\omega(G) = |K^n| + \ldots + |K^1| + \omega(G_0)$ and $\chi(G) = |K^n| + \ldots + |K^1| + \chi(G_0)$. 
\end{prop}
\begin{proof}
    We prove the base case $n=1$; the full result holds by induction. 
    
    Let $A$ be a maximum clique in $G_0$. The set $K \cup A$ induces a clique in $G$, so $\omega(G) \ge |K| + \omega(G_0)$. Any maximum clique including $s \in S$ contains at most $|K| + 1$ vertices, which is not greater than $|K| + \omega(G_0)$, and any maximum clique including $v \in V(G_0) - A$ still contains at most $|K| + \omega(G_0)$ vertices. Hence $\omega(G) = |K| + \omega(G_0)$. 

    A proper coloring of $G_0$ must use at least $\chi(G_0)$ colors, and a proper coloring of $G$ must give every vertex in $K$ a color distinct from one another and from the colors used on $G_0$. Hence $\chi(G) \ge |K| + \chi(G_0)$. All vertices in $S$ can be colored with any color used on $G_0$, since $G_0$ is non-empty. This produces a proper coloring of $G$ using $|K| + \chi(G_0)$ colors, so we conclude that $\chi(G) = |K| + \chi(G_0)$.
\end{proof}

Where $G_n, \ldots, G_1$ each have exactly one split partition, $\omega(G) = \omega(G_n) + \ldots + \omega(G_0)$ and $\chi(G) = \chi(G_n) + \ldots + \chi(G_0)$. 

Recall that a graph is apex-perfect if and only if there exists $v \in V(G)$ such that the induced subgraph $G - \{v\}$ is perfect. 

\begin{theorem}
    \label{prop:chi_bound} 
    The hereditary closure of the unigraphs is a subset of the apex-perfect graphs; hence the class is $\chi$-bounded by the function $f(x) = x+1$.  
\end{theorem}
\begin{proof}
    All indecomposable unigraphs are perfect except for $C_5$, which is apex-perfect. Let $G$ be a unigraph, and write $G$ as the composition $G_n \circ \ldots \circ G_1 \circ G_0$ for some set of indecomposable unigraphs $G_0, \ldots, G_n$. Per \cite{Ty00}, all indecomposable split unigraphs have exactly one split partition, so the composition is well-defined. As split graphs, all of $G_1, \ldots, G_n$ are perfect. The graph $G_0$ is either perfect or isomorphic to $C_5$, such that $\chi(G_0) \le \omega(G_0) + 1$. Thus $\chi(G) = \chi(G_n) + \ldots + \chi(G_0) \le \omega(G_n) + \ldots + \omega(G_1) + \omega(G_0) + 1$, by Proposition \ref{prop:composition_coloring}.  
    
    Since the set of apex-perfect graphs is hereditary, we conclude that any graph induced in a unigraph is apex-perfect. Hence the hereditary closure of the unigraphs is a subset of the apex-perfect graphs.
\end{proof}

\section{Hereditary Perfect- and Chordal-Unigraphs}

While all unigraphs are apex-perfect, not all are perfect ($C_5$ is the most notable counterexample). In this section we attempt to identify ``better" related classes of graphs by considering hereditary perfect-unigraphs and hereditary chordal-unigraphs. The characterization of the former class is straightforward. 

\begin{theorem}
    \label{prop:hereditary_perfect_graphs}
    A graph $G$ is a hereditary perfect-unigraph if and only if $G$ contains none of the graphs in $\mathcal{F}_{H(\mathcal{U})}$ nor $C_5$ as an induced subgraph. 
\end{theorem}
\begin{proof}
    Suppose $G$ is a hereditary perfect-unigraph. Of course it cannot contain $C_5$; moreover, all graphs in $\mathcal{F}_{H(\mathcal{U})}$ are perfect but not hereditary unigraphs, so $G$ contains none of them as an induced subgraph. 

    For the reverse, suppose $G$ contains no element of $\mathcal{F}_{H(\mathcal{U})}$ nor $C_5$ as an induced subgraph. Thus $G$ is a hereditary unigraph, and since $\{P_5, \overline{P_5}\} \subset \mathcal{F}_{H(\mathcal{U})}$, it follows that $G$ is perfect. Since $\mathcal{F}_{H(\mathcal{U})}$ is closed under degree sequence realizations, it follows that $G$ is a hereditary perfect-unigraph. 
\end{proof}

\begin{figure}
\centering
    \includegraphics[height=3cm]{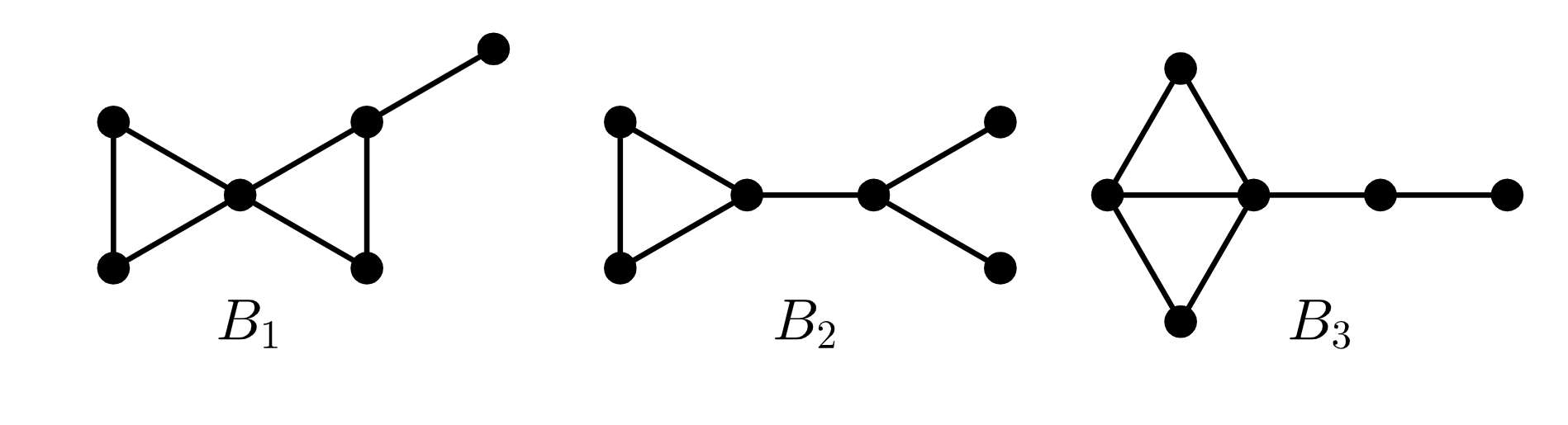}
  \caption{Three of the thirteen forbidden induced subgraphs of the hereditary chordal-unigraphs.}
\label{fig:chordal_forbidden}
\end{figure}

We provide two standard theorems about split graphs in advance of their use in the following theorem characterizing $H(\mathcal{U}_C)$. 

\begin{theorem} \cite{FoHa77}
\label{prop:split_mfis}
    A graph $G$ is a split graph if and only if it contains none of $C_4$, $C_5$, or $2K_2$ as an induced subgraph. 
\end{theorem}

\begin{theorem} \cite{HaSi81}
\label{prop:split_degree}
    A graph $G$ is a split graph if and only if all other realizations of its degree sequence are split graphs.
\end{theorem}

Let $R, S$ be as defined before Theorem \ref{prop:hereditary_unigraphs}. Let $B_1, B_2, B_3$ be as shown in Figure \ref{fig:chordal_forbidden}. Let $\mathcal{F}_{H(\mathcal{U_C})} = \{C_4, C_5, P_5, K_3 + K_2, 2P_3, P_4 + K_2, R, \overline{R}, S, \overline{S}, B_1, B_2, B_3\}$. 

\begin{theorem}
    \label{prop:hereditary_chordal_graphs}
    A graph $G$ is a hereditary chordal-unigraph if and only if $G$ contains none of $\mathcal{F}_{H(\mathcal{U_C})}$ as an induced subgraph. 
\end{theorem}
\begin{proof}
    Suppose $G$ is a hereditary chordal-unigraph. Since $G$ is chordal, it contains neither $C_4$ nor $C_5$ as an induced subgraph. The graphs $P_5$, $K_3 + K_2$, $2P_3$, $P_4 + K_2$, $R$, $\overline{R}$, $S$, and $\overline{S}$ are all chordal and minimal non-unigraphs, so $G$ can contain none of these. The graphs $B_1$ and $B_2$ are both chordal realizations of $(4,3,2,2,2,1)$. The sequence $(3,3,2,2,1,1)$ has two chordal realizations - one is $B_3$ and the other contains $P_5$ as an induced subgraph - so $G$ cannot contain $B_3$ as an induced subgraph. 

    Suppose $G$ contains none of $\mathcal{F}_{H(\mathcal{U_C})}$ as an induced subgraph. Since $G$ cannot contain $C_4, C_5, $ or $P_5$, $G$ is chordal. If $G$ is a split graph, then all realizations of its degree sequence are split by Theorem \ref{prop:split_degree}, and $G$ is a hereditary chordal-unigraph if and only if $G$ is a hereditary unigraph. Since $G$ contains none of $R$, $\overline{R}$, $S$, or $\overline{S}$, it follows that $G$ is a split hereditary unigraph by Theorem \ref{prop:hereditary_unigraphs}. If $G$ is a chordal-unigraph and contains an isolated or dominating vertex $v$, then $G - \{v\}$ is also a chordal-unigraph, since $v$ is isolated or dominating in every realization, and cannot be induced in a cycle. Hence we assume that $G$ contains neither isolated nor dominating vertices. 

    Hence we assume $G$ is not a split graph, and therefore must contain a copy of $2K_2$ as an induced subgraph by Theorem \ref{prop:split_mfis}. Suppose vertices $a,b,c,d$ of $G$ induce $2K_2$ with edges $ab$ and $cd$. If another vertex $v \in V(G)$ is adjacent, up to symmetry, to $a$ and $b$ but not $c$ or $d$, then $\{a,b,c,d,v\}$ induces $K_3+K_2$, for a contradiction. If another vertex $v \in V(G)$ is adjacent, up to symmetry, to $a$ and $c$ but not $b$ or $d$, then $\{a,b,c,d,v\}$ induces $P_5$, for a contradiction. Hence no vertex in $G$ is adjacent to exactly two vertices out of $a,b,c,$ and $d$. 
    
    If other vertices $v_1, \ldots, v_k \in V(G)$ are adjacent to exactly one vertex from $\{a,b,c,d\}$, the resulting graph contains $P_5$, $P_4+K_2,$ or $2P_3$ unless $v_1, \ldots, v_k$ are all adjacent to the same vertex (say, $a$) and not adjacent to one another. The resulting graph is isomorphic to the disjoint union of a star and $K_2$, and is unigraphic by Theorem \ref{prop:indecomposable_unigraphs}.

    If other vertices $v_1, \ldots, v_k \in V(G)$ are adjacent to exactly three vertices from $\{a,b,c,d\}$, then $v_1, \ldots, v_k$ are all adjacent to one another and each is adjacent to $a, d$, and either $b$ or $c$, up to symmetry; otherwise, $G$ contains $C_4$ or $C_5$ as an induced subgraph. 

    Suppose $v_1, \ldots, v_k$ are adjacent to $a, b, $ and $d$, and vertices $w_1, \ldots, w_m$ are adjacent to $a, c$, and $d$. The degree sequence of $G$ is thus $((k+m+2)^{k+m}, (k+m+1)^2, k+1, m+1)$. We claim no other chordal realization $G'$ of this sequence exists. First, such a realization cannot contain $C_4$ or $C_5$. Since $G$ is not a split graph, $G'$ cannot be split, so it must contain $2K_2$ as an induced subgraph. Since vertices of degree $k+m+2$ are adjacent to all but one other vertex in $G'$, none can be part of an induced $2K_2$, so the four lower-degree vertices must form this $2K_2$. Second, if two vertices $p,q$ of degree $k+m+2$ are not adjacent to one another, then they must be complete to the remaining vertices. Together with any two non-adjacent vertices $r,s$, the set $\{p,q,r,s\}$ induces $C_4$ in $G'$, for a contradiction. Hence the vertices of degree $k+m+2$ are complete to one another, and each is non-adjacent to exactly one more vertex. Third, if $k\ge 0$, $m \ge 0$, and the vertex of degree $k+1$ is adjacent to the vertex of degree $m+1$, then each is non-adjacent to some vertex of degree $k+m+2$, which must be distinct. These four vertices then induce $C_4$, for a contradiction. If $k=0$ or $m=0$, then there is only one possible induced $2K_2$, up to symmetry. Hence one vertex of degree $k+m+1$ is adjacent to the vertex of degree $k+1$, and the other is adjacent to the vertex of degree $m+1$, inducing $2K_2$. The vertices of degree $k+m+1$ are by necessity complete to the vertices of degree $k+m+2$, each of which is adjacent to either the vertex of degree $k+1$ or the vertex of degree $m+1$, but not both, producing $G$. 

    If there exist $v$ adjacent to one of $\{a,b,c,d\}$ and $w$ adjacent to three of $\{a,b,c,d\}$, then $G$ must contain one of $P_5$, $C_4$, $B_1$, $B_2$, or $B_3$ as an induced subgraph. Hence $G$ cannot contain vertices with both one and three edges to vertices in $\{a,b,c,d\}$, aside from $a,b,c$, and $d$. 

    If there exists a vertex $v \in V(G)$ adjacent to all of $\{a,b,c,d\}$, then since $G$ has no dominating vertices, there exists $w \in V(G)$ not adjacent to $v$. If $w$ is adjacent to one of $\{a,b,c,d\}$, then $G$ contains $B_1$ as an induced subgraph. If $w$ is adjacent to three or four of $\{a,b,c,d\}$, then $G$ contains $C_4$ as an induced subgraph. Hence $w$ is adjacent to none of $\{a,b,c,d\}$, and so there exists $x$ adjacent to $w$ and none of $\{a,b,c,d,v\}$. If $x$ has any other neighbor in $G$, then $G$ contains $2P_3$ or $K_3 + K_2$ as an induced subgraph, so instead $\{w,x\}$ is a disconnected component of $G$.  Thus $G$ contains $K_3 + K_2$ as an induced subgraph, for a contradiction. 

    If there exists a vertex $v \in V(G)$ adjacent to none of $\{a,b,c,d\}$, then since $G$ has no isolated vertices, there exists $w \in V(G)$ adjacent to $v$. If $w$ is adjacent to one of $\{a,b,c,d\}$, then $G$ contains $P_4+K_2$ as an induced subgraph. If $w$ is adjacent to all of $\{a,b,c,d\}$, we use the argument in the above paragraph to show $G$ contains $K_3 + K_2$ as an induced subgraph. If $w$ is adjacent to none of $\{a,b,c,d\}$, then $G$ is isomorphic to $nK_2$ and is a hereditary unigraph by Theorem \ref{prop:indecomposable_unigraphs}.
    
    Suppose instead, then, that $w$ is adjacent to three of $\{a,b,c,d\}$, then it must be the only vertex adjacent to three of $\{a,b,c,d\}$, else $G$ contains $R$ or $\overline{R}$ as an induced subgraph. If more than one vertex of $G$ is not adjacent to $w$, then $G$ contains either $B_3$ or $P_5$ as an induced subgraph. The graph $G - \{w\}$ contains several disconnected components; since one of them, $\{a,b\}$, induces $K_2$, the others cannot contain an induced $P_4$ or $K_3$. Hence each is isomorphic to $K_1, K_2,$ or $P_3$, and at most one is isomorphic to $P_3$. Should the one vertex not adjacent to $w$ be in a component isomorphic to $P_3$, $G$ contains either $B_1$ or $C_4$ as an induced subgraph. The result is a unigraph: $w$ is adjacent to all vertices but one of degree $1$; there is at most one vertex of degree $3$, which forms an induced $P_3$ in $G - \{w\}$ together with two vertices of degree $2$, and all remaining vertices of degree $1$ and $2$ induce the unigraph $xK_2 + yK_1$ in $G - \{w\}$, for $x \ge 2$ and $y \ge 0$.     
\end{proof}

\end{document}